\documentclass[12pt]{article}
\usepackage[utf8]{inputenc}
\usepackage{color}
\usepackage{amsmath,amssymb,amsthm,amsfonts}
\usepackage{float}
\usepackage{graphicx}
\usepackage{comment}
\usepackage{multirow}
\RequirePackage{hyperref}

\newtheorem{remark}{Remark}
\newtheorem{corollary}{Corollary}
\newtheorem{definition}{Definition}
\newtheorem{theorem}{Theorem}

\usepackage{geometry}
\begin{document}
\begin{center}
\textsc{{ \Large A fractional order cubic differential equation. Applications to natural resource management}\\[0pt]}
\vspace{0.4cm} \hspace{0.2cm} Melani Barrios$^{1,2}$ \hspace{0.5cm} Gabriela Reyero$^{1}$ \hspace{0.5cm} Mabel Tidball$^{3}$ \\[0pt]
\end{center}

\scriptsize                                                                               
 $\,^{1}$  Departamento de Matem\'atica, Facultad de Ciencias Exactas, Ingeniería y Agrimensura, Universidad Nacional de Rosario, Avda. Pellegrini $250$, S$2000$BTP Rosario, Argentina.
 
$\,^{2}$ CONICET, Departamento de Matem\'atica, Facultad de Ciencias Exactas, Ingeniería y Agrimensura, Universidad Nacional de Rosario, Avda. Pellegrini $250$, S$2000$BTP Rosario, Argentina.

$\,^{3}$ CEE-M, Universidad de Montpellier, CNRS, INRA, SupAgro, Montpellier, France.\\

\normalsize

Correspondence should be addressed to melani@fceia.unr.edu.ar\\

\begin{abstract} 
 An analysis of a fractional cubic differential equation is presented, which is a generalization of different versions of fractional logistic equations, in order to obtain simpler numerical methods that globalize and extend the results already obtained, and allow the comparison between several methods. The existence and uniqueness of the solution of the fractional cubic problem subject to an initial value are demonstrated, a stability analysis is performed and, based on the implementation of numerical methods, a comparison is made between the different logistic models.
\end{abstract}


\textbf{Keywords} Fractional derivatives and integrals, Fractional differential equations, Environmental economics, Stability.\\

\textbf{Mathematics Subject Classification} 26A33, 34A08, 34D20, 91B76.

\section{Introduction}

Predicting the future of a population number is one of the most important factors needed for the good management of it. This has been treated by several known methods, one of them being the development of a mathematical model which describes the population growth. The model generally takes the form of a differential equation, or a system of differential equations, according to the complexity of the underlying properties of the population. The currently most used growth models are those that have a sigmoid solution in time, including Gompertz and Verhulst's logistic equations. The logistic equation is commonly used in population growth models, disease propagation epidemic, and social networks \cite{AmNuAnSu, Cl}.

Motivated by its applications in different scientific areas (electricity, magnetism, mechanics, fluid dynamics, medicine, etc. \cite{Alm, BaRe2, GoRe, Goo, Hil, Kil}), fractional calculus is in development, which has led to great growth in its study in recent decades. The fractional derivative is a nonlocal operator \cite{Die, Pod}, this makes fractional differential equations good candidates for modeling situations in which it is important to consider the history of the phenomenon studied \cite{FeSa}, unlike the models with classical derivative where this is not taken into account. There are several definitions of fractional derivatives. The most commonly used are the Riemann-Liouville fractional derivative and the Caputo fractional derivative. It is important to note that while the Riemann-Liouville fractional derivative \cite{Old} is historically the most studied approach to fractional calculus, the Caputo fractional derivative is more popular among physicists and scientists because the formulation of initial values problems with this type of derivative is more similar to the formulation with classical derivative.

There have been different fractional versions of the logistic equation, among which we find: the usual fractional logistic equation \cite{EsEmEs, KSQB}, and the fractional logistic equation with the Allee effect \cite{SyMaSh}. In this paper, the usual fractional logistic equation and with Allee effect will be analyzed, considering in both the harvest of the studied resource. For this reason, a study of a fractional cubic equation will be presented, which is a generalization of these different versions of fractional logistic equations. In the first part of the work, existence and uniqueness of the solution of the fractional cubic problem subject to an initial value will be proved, in the second part a stability analysis will be performed and finally, from the implementation of numerical methods, a comparison between the different logistic models will be seen.

\section{Preliminaries}\label{sec2}
\subsection{Different types of fractional logistic equations and fractional cubic equation}

\begin{itemize}
\item \textbf{Fractional logistic equation:}
\begin{equation}
\,_0^{C} D_t^{\alpha}\left[x\right](t)=rx(t)\left(1-\frac{x(t)}{K}\right),
\label{eclogistica}
\end{equation}
where $r>0$ is the intrinsic growth rate and $K>0$ the carrying capacity of the resource.
\item \textbf{Fractional logistic equation with harvest:}
\begin{equation}
\,_0^{C} D_t^{\alpha}\left[x\right](t)=rx(t)\left(1-\frac{x(t)}{K}\right)-Ex(t),
\label{eclogisticacosecha}
\end{equation}
where $E>0$ is the effort to be made to extract a proportion of the resource studied.
\item \textbf{Fractional logistic equation with Allee efect:}
\begin{equation}
\,_0^{C} D_t^{\alpha}\left[x\right](t)=rx(t)\left(1-\frac{x(t)}{K}\right)(x(t)-m),
\label{eclogisticaallee}
\end{equation}
where $m>0$ is the threshold of the Allee effect, which means the minimum population density for the growth of certain species, below which the population dies out (the population growth rate is positive only within the range $ m <x <K $ and is negative outside this range).
\item \textbf{Fractional logistic equation with Allee efect with harvest:}
\begin{equation}
\,_0^{C} D_t^{\alpha}\left[x\right](t)=rx(t)\left(1-\frac{x(t)}{K}\right)(x(t)-m)-Ex(t).
\label{eclogisticaalleecos}
\end{equation}
\end{itemize}

The objective of this paper is to study a fractional cubic equation, which is a generalization of these different versions of fractional logistic equations. Then, consider the fractional cubic equation given by:
\begin{equation}
\,_0^{C} D_t^{\alpha}\left[x\right](t)=ax^3(t)+bx^2(t)+cx(t),
\label{cubicafrac}
\end{equation}
where $a,b,c \in \mathbb{R}$.

It is necessary to introduce the meaning of fractional derivative and some of its properties.

\subsection{Fractional Calculus} \label{sec:frac}
\begin{definition}
The Gamma function, $\Gamma: (0, \infty)\rightarrow \mathbb{R}$, is defined by:
\begin{equation}
 \Gamma(t) = \int_{0}^{\infty} s^{t-1} e^{-s} \, ds.
\label{gamma}
\end{equation}
\end{definition}
\begin{definition}
The Riemann-Liouville fractional integral operator of order $\alpha \in \mathbb{R}^{+}_{0}$ is defined in $L^1[a,b]$ by:
\begin{equation}
 \,_{a}I_{t}^{\alpha} [f] (t) = \dfrac{1}{\Gamma(\alpha)} \int_{a}^{t} (t-s)^{\alpha -1} f(s) \, ds.
\label{frac1}
\end{equation}
\end{definition}
\begin{definition}\label{defi}
The Caputo fractional derivative operator of order $\alpha \in \mathbb{R}^{+}_{0}$ is defined by:
\begin{equation}
 \,_{a}^{C}D_{t}^{\alpha}[f] (t)=\left(\,_{a}I_{t}^{n-\alpha} \circ \dfrac{d^{n}}{dt^{n}}\right)[f] (t)
\label{frac5}
\end{equation}
as long as $\dfrac{d^{n}f}{dt^{n}} \in L^1[a,b]$, and $n=\left\lceil \alpha \right\rceil$.
\end{definition}

In \cite{OdKuShAlHa} the following result is proved:

\begin{theorem} [Generalized Mean Value Theorem] \label{thm:tvmg} 
Let $f(t)\in AC[0,b]$ (the set of absolutely continuous functions on $[0,b]$). Then, for $0<\alpha\leq 1$:
\[f(t)=f(0)+\frac{1}{\Gamma(\alpha+1)} \,_{0}^{C}D_{t}^{\alpha}[f](\xi)t^{\alpha},\]
\noindent with $0\leq\xi\leq t, \; \forall t\in[0,b]$.
\end{theorem}

\begin{remark}
When $\alpha=1$, the generalized mean value theorem reduces to the classical mean value theorem.
\end{remark}

\begin{corollary} \label{cor:crec}
Suppose that $f(t) \in AC[0,b]$ and $\,_{0}^{C}D_{t}^{\alpha}[f](t) \in C(0,b]$ for $0<\alpha\leq 1$. If  $\,_{0}^{C}D_{t}^{\alpha}[f](t)\geq 0 \, \left(\,_{0}^{C}D_{t}^{\alpha}[f](t)> 0\right)$, $\forall t \in (0,b)$, then $f(t)$ is non-decreasing (increasing) and if $\,_{0}^{C}D_{t}^{\alpha}[f](t)\leq 0 \, \left(\,_{0}^{C}D_{t}^{\alpha}[f](t)< 0\right)$, $\forall t \in (0,b)$, then $f(t)$ is non-increasing (decreasing) for all $t \in [0,b]$.
\end{corollary}

\section{Study of the fractional cubic equation}
\subsection{Existence and uniqueness of the solution}\label{sec:uniq}
\,\\
Consider the following initial value problem
\begin{equation}
 \left\{ \begin{array}{l}
             \,_0^{C} D_t^{\alpha}\left[x\right](t)=ax^3(t)+bx^2(t)+cx(t), \\
             \\ x(0)=x_0, \\
             \end{array}
   \right.
	\label{problemappal}
\end{equation}
with $t>0$ and $0<\alpha\leq 1$.
\begin{definition}
The Bielecki norm is defined by:
\[\left\|x\right\|_N=\sup_{t \in I} \left|e^{-Nt}x(t)\right|\]
with $N>0$ and $I=\left[0,T\right]$, \cite{Me}.
\end{definition} 
\begin{definition} $x(t)$ is solution of the initial value problem (\ref{problemappal}) if:
\begin{enumerate}
\item $(t,x(t)) \in I\times H$ with $I=\left[0,T\right]$ and $H=\left[-h,h\right]$.
\item $x(t)$ satisfies (\ref{problemappal}).
\end{enumerate}
\end{definition}
\begin{theorem}
If $S=\left\{x \in L^1\left[0,T\right], \, \left\|x\right\|=\left\|e^{-Nt}x(t)\right\|_{L^1}\right\}$ and $N>0$ is such that $N^{\alpha}>3\left|a\right|h^2+2\left|b\right|h+\left|c\right|$ then the problem  (\ref{problemappal}) has unique solution $x \in C(I)$, with $x' \in S$. 
\end{theorem}  
\begin{proof}
By definition \ref{defi}, 
\[ \,_0^{C} D_t^{\alpha}\left[x\right](t)=\,_0 I_t^{1-\alpha} \frac{d}{dt}\left[x\right](t)=ax^3(t)+bx^2(t)+cx(t).\]
Applying $\,_0 I_t^{\alpha}$ in both members, it is obtained
\begin{equation}
x(t)=x_0+\,_0 I_t^{\alpha} \left[ax^3+bx^2+cx\right](t).
\label{imp1}
\end{equation}
Let $G:C(I)\rightarrow C(I)$ be an operator such as $G(x(t))= x_0+\,_0 I_t^{\alpha} \left[ax^3+bx^2+cx\right](t)$. \\
To see that $G$ has only one fixed point, it can be proved that $G$ is contractive with the Bielecki norm $\left\|.\right\|_N$.
\[e^{-Nt}(Gx-Gy)=e^{-Nt} \,_0 I_t^{\alpha} \left[a(x^3-y^3)+b(x^2-y^2)+c(x-y)\right]=\]
\[=\int_0^t \frac{(t-s)^{\alpha-1}}{\Gamma(\alpha)} e^{-Nt} \left[a(x^3(s)-y^3(s))+b(x^2(s)-y^2(s))+c(x(s)-y(s))\right]ds=\]
\[\displaystyle\int_0^t \dfrac{(t-s)^{\alpha-1}}{\Gamma(\alpha)} e^{-Nt+Ns} e^{-Ns} \left[a(x^3(s)-y^3(s))+b(x^2(s)-y^2(s))+c(x(s)-y(s))\right]ds=\]

 \[\displaystyle\int_0^t \frac{(t-s)^{\alpha-1}}{\Gamma(\alpha)} e^{-N(t-s)} \left[e^{-Ns}(x(s)-y(s))\right] \left[a(x^2(s)+x(s)y(s)+y^2(s))+b(x(s)+y(s))+c\right]ds\leq\]

\[\leq \left\|x-y\right\|_{N} \displaystyle\int_0^t  \frac{(t-s)^{\alpha-1}}{\Gamma(\alpha)} e^{-N(t-s)} \left[a(x^2(s)+x(s)y(s)+y^2(s))+ +b(x(s)+y(s))+c\right]ds\leq \]


 with the change of variable $u=t-s$, it results
 \[\left\|x-y\right\|_{N} \int_0^t \frac{u^{\alpha-1}}{\Gamma(\alpha)} e^{-Nu}  \left[a(x^2(t-u)+x(t-u)y(t-u)+y^2(t-u))+b(x(t-u)+y(t-u))+c\right]du.\]
 
 Applying module in both members,
\[\left|e^{-Nt}(Gx-Gy)\right|\leq \left\|x-y\right\|_N \left(3\left|a\right|h^2+2\left|b\right|h+\left|c\right|\right)
N^{-\alpha}\int_0^t N^{\alpha}\frac{(u)^{\alpha-1}}{\Gamma(\alpha)} e^{-Nu}du.\]

By Gamma function definition, it is obtained
 
\[\left|e^{-Nt}(Gx-Gy)\right|\leq \left\|x-y\right\|_N \frac{3\left|a\right|h^2+2\left|b\right|h+\left|c\right|}{N^{\alpha}}
\frac{\Gamma(\alpha)}{\Gamma(\alpha)},\]

then

\[\left|e^{-Nt}(Gx-Gy)\right|\leq \left\|x-y\right\|_N \frac{3\left|a\right|h^2+2\left|b\right|h+\left|c\right|}{N^{\alpha}}.\]

 Considering $N$ such as $N^{\alpha}>3\left|a\right|h^2+2\left|b\right|h+\left|c\right|$, then  $\left\|Gx-Gy\right\|_N< \left\|x-y\right\|_N$, therefore the operator $G$ results contractive and then by the fixed point theorem, the problem has unique solution.\\
 To see that $x \in C(I)$ and $x' \in S$, from equation (\ref{imp1}),
\[x(t)=x_0+\,_0 I_t^{\alpha} \left[ax^3+bx^2+cx\right](t),\]

and it can be written as

\[x(t)=x_0+\,_0 I_t^{\alpha} \left[\,_0 I_t^{1}\frac{d}{dt}\left[ax^3+bx^2+cx\right](t)+ax_0^3+bx_0^2+cx_0\right].\]

Applying linearity and the Riemann-Liouville integral definition, it is obtained

\[x(t)=x_0+\frac{t^{\alpha}}{\Gamma(\alpha+1)}\left(ax_0^3+bx_0^2+cx_0\right)+\,_0 I_t^{\alpha+1}
\left[3a(x)^2x'+2bxx'+cx'\right](t),\]

then $x \in C(I)$. Deriving respect to $t$,

\[x'(t)=\frac{t^{\alpha-1}}{\Gamma(\alpha)}\left(ax_0^3+bx_0^2+cx_0\right)+\,_0 I_t^{\alpha}
\left[3a(x)^2x'+2bxx'+cx'\right](t),\]

which means that,

\[e^{-Nt}x'(t)=e^{-Nt}\left(\frac{t^{\alpha-1}}{\Gamma(\alpha)}\left(ax_0^3+bx_0^2+cx_0\right)+\,_0 I_t^{\alpha}
\left[3a(x)^2x'+2bxx'+cx'\right](t)\right),\]

where it can be deduced that $e^{-Nt}x'(t) \in L^{1}[0,T]$ and therefore $x' \in S$.

To prove that the expression of $x(t)$ given by (\ref{imp1}) verify the problem (\ref{problemappal}), deriving respect to $t$,

\[\frac{d}{dt}x(t)=\dfrac{d}{dt}\,_0 I_t^{\alpha} \left[ax^3+bx^2+cx\right](t).\]

Applying $\,_0 I_t^{1-\alpha}$,

\[\,_0 I_t^{1-\alpha} \frac{d}{dt}x(t)=\,_0 I_t^{1-\alpha}\frac{d}{dt}\,_0 I_t^{\alpha}\left[ax^3+bx^2+cx\right](t).\]

By Caputo derivative in definition \ref{defi},

\[\,_0^C D_t^{\alpha} x(t)= \frac{d}{dt}\,_0 I_t^{1-\alpha}\,_0 I_t^{\alpha} \left[ax^3+bx^2+cx\right](t),\]

\[\,_0^C D_t^{\alpha} x(t)= \frac{d}{dt}\,_0I_t \left[ax^3+bx^2+cx\right](t),\]

and then, it is the fractional cubic equation

\[\,_0^C D_t^{\alpha} x(t)=\left[ax^3+bx^2+cx\right](t).\]

Finally, to verify the initial condition, from the expression in (\ref{imp1}),

\[x(0)=x_0+\underbrace{\,_0 I_0^{\alpha} \left[ax^3+bx^2+cx\right](t)}_{=0}\]

then, 
\[x(0)=x_0.\]

Accordingly, integral equation (\ref{imp1}) is equivalent to problem (\ref{problemappal}) and the theorem is proved.
\end{proof}

\subsubsection{Existence of solution for fractional logistic equations}
\,\\
Following, it will be shown that the condition for $N$, $N^{\alpha}>3\left|a\right|h^2+2\left|b\right|h+\left|c\right|$, are the same as the conditions for the existence of fractional logistic equations.  
\begin{itemize}
\item \textbf{Fractional logistic equation:}
\begin{equation}
\,_0^{C} D_t^{\alpha}\left[x\right](t)=rx(t)\left(1-\frac{x(t)}{K}\right).
\label{eclogistica}
\end{equation}
This is equation (\ref{problemappal}) with $a=0$, $b=-\tfrac{r}{K}$ and $c=r$. 

Then $N^{\alpha}>2\frac{r}{K}h+r=r\left(\frac{2h}{K}+1\right)$. This result is analogous to that obtained in \cite{EsEmEs}. 
\item \textbf{Fractional logistic equation with harvest:}
\begin{equation}
\,_0^{C} D_t^{\alpha}\left[x\right](t)=rx(t)\left(1-\frac{x(t)}{K}\right)-Ex(t).
\label{eclogisticacosecha}
\end{equation}
This is equation (\ref{problemappal}) with $a=0$, $b=-\tfrac{r}{K}$ and $c=r-E$.

Then $N^{\alpha}>2\frac{r}{K}h+(r-E)$.
\item \textbf{Fractional logistic equation with Allee effect:}
\begin{equation}
\,_0^{C} D_t^{\alpha}\left[x\right](t)=rx(t)\left(1-\frac{x(t)}{K}\right)(x(t)-m).
\label{eclogisticaallee}
\end{equation}
This is equation (\ref{problemappal}) with $a=-\tfrac{r}{K}$, $b=\left(\tfrac{m}{K}+1\right)r$ and $c=-rm$. 

Then $N^{\alpha}>3\frac{r}{K}h^2+2\left(\tfrac{m}{K}+1\right)r h+rm=r\left(m+\left(\tfrac{m}{K}+1\right)2h+3\frac{h^2}{K}\right)$.
\item \textbf{Fractional logistic equation with Allee effect with harvest:}
\begin{equation}
\,_0^{C} D_t^{\alpha}\left[x\right](t)=rx(t)\left(1-\frac{x(t)}{K}\right)(x(t)-m)-Ex(t).
\label{eclogisticaalleecos}
\end{equation}
This is equation (\ref{problemappal}) with $a=-\tfrac{r}{K}$, $b=\left(\tfrac{m}{K}+1\right)r$ and $c=-rm-E$.

Then $N^{\alpha}>3\frac{r}{K}h^2+2\left(\tfrac{m}{K}+1\right)r h+rm+E=r\left(m+\left(\tfrac{m}{K}+1\right)2h+3\frac{h^2}{K}\right)+E$.
\end{itemize}

\subsection{Stability analysis of the fractional cubic equation}
\subsubsection{General case}
\,\\
Consider the following fractional initial value problem,
\begin{equation}
 \left\{ \begin{array}{l}
             \,_0^{C} D_t^{\alpha}\left[x\right](t)=f(x(t)), \\
             \\ x(0)=x_0, \\
             \end{array}
   \right.
	\label{cfgeneral}
\end{equation}

with $t>0$ and $0<\alpha\leq 1$.

To find the equilibrium points $x_{eq}$ of this equation, the condition $\,_0^{C} D_t^{\alpha}x(t)=0$ is stated, which means,
\begin{equation}
f(x_{eq})=0.
\label{geq}
\end{equation}

To study the stability of each equilibrium point, consider the jacobian matrices $A_{eq}$ of $f$ evaluated at the equilibrium points,
\[A_{eq}=\left[\left.f'(x(t))\right|_{x_{eq}}\right].\]

The eigenvalues $\lambda_{eq}$ of $A_{eq}$, are calculated for each equilibrium point.

The following theorem can be seen in \cite{AhESES} and \cite{Die}.
\begin{theorem}
Let $\lambda_{eq}$ be the eigenvalues of $A_{eq}$, jacobian matrix associated to each equilibrium point $x_{eq}$. Then
\begin{itemize}
\item if $arg(\lambda_{eq})< \frac{\alpha\pi}{2}$, then the equilibrium point $x_{eq}$ is locally unstable (U),
\item if $arg(\lambda_{eq})\geq \frac{\alpha\pi}{2}$, then the equilibrium point $x_{eq}$ is locally stable (S), being locally asymptotically stable (AS) if ${arg(\lambda_{eq})> \frac{\alpha\pi}{2}}$.
\end{itemize}
\end{theorem}

 Now, the equilibrium points and their stability of the fractional cubic equation will be examined.
 
\subsubsection{Fractional cubic equation case}
Consider $0<\alpha \leq 1$ and $a,b,c \in \Bbb{R}$ for the following problem
\begin{equation}
 \left\{ \begin{array}{l}
             \,_0^{C} D_t^{\alpha}\left[x\right](t)=f(x(t))=ax^3(t)+bx^2(t)+cx(t), \\
             \\ x(0)=x_0. \\
             \end{array}
   \right.
	\label{problemappal2}
\end{equation}

In the Table \ref{table:pointsofequilibrium}, the different points of equilibrium and the stabilities of each one can be observed. To see the complete analysis refer to appendix or \cite{BaReTid}.

\begin{table}[!h]
\centering
\caption{Points of equilibrium \label{table:pointsofequilibrium}}
\begin{tabular}{|l|l|l|l|l|}
\hline
\multirow{10}{*}{$a\neq 0$} & \multirow{2}{*}{$x_1 = 0$} & \multicolumn{2}{c|}{$c>0$ } & $x_1$ is (U)\\
& & \multicolumn{2}{c|}{$c<0$ } & $x_1$ is (AS)\\ \cline{2-5}
& \multirow{4}{*}{$x_2=\frac{-b+\sqrt{b^2-4ac}}{2a}$} & \multirow{2}{*}{$a>0$} & $b<\sqrt{b^2-4ac}$ & $x_2$ is (U)\\
& & & $b>\sqrt{b^2-4ac}$ & $x_2$ is (AS)\\ \cline{3-5}
& & \multirow{2}{*}{$a<0$} & $b<\sqrt{b^2-4ac}$ & $x_2$ is (AS)\\
& & & $b>\sqrt{b^2-4ac}$ & $x_2$ is (U)\\ \cline{2-5}
& \multirow{4}{*}{$x_3=\frac{-b-\sqrt{b^2-4ac}}{2a}$} & \multirow{2}{*}{$a>0$} & $b<-\sqrt{b^2-4ac}$ & $x_3$ is (AS)\\
& & & $b>-\sqrt{b^2-4ac}$ & $x_3$ is (U)\\ \cline{3-5}
& & \multirow{2}{*}{$a<0$} & $b<-\sqrt{b^2-4ac}$ & $x_3$ is (U)\\
& & & $b>-\sqrt{b^2-4ac}$ & $x_3$ is (AS)\\ \hline
\hline
\multirow{4}{*}{$a = 0$} & \multirow{2}{*}{$x_1 = 0$} & \multicolumn{2}{c|}{$c>0$ } & $x_1$ is (U)\\
& & \multicolumn{2}{c|}{$c<0$ } & $x_1$ is (AS)\\ \cline{2-5}
& \multirow{2}{*}{$x_2 = -\frac{c}{b}$} & \multicolumn{2}{c|}{$c>0$ } & $x_2$ is (AS)\\
& & \multicolumn{2}{c|}{$c<0$ } & $x_2$ is (U)\\ \hline

\end{tabular}
\end{table}

\subsubsection{Equilibrium and stability applied to fractional logistic equations} 
\,\\
Following, it will be shown that the results found, in the previous subsection, coincide with the results of the fractional logistic equations (\ref{eclogistica}), (\ref{eclogisticacosecha}), (\ref{eclogisticaallee}) and (\ref{eclogisticaalleecos}). 
\begin{itemize}
\item \textbf{Fractional logistic equation:}
This is equation (\ref{problemappal}) with $a=0$, $b=-\tfrac{r}{K}$ and $c=r$.
Then, $a=0$, $b<0$ and $c>0$, therefore the equilibrium points are $x_1=0$ which is (U) and  $x_2=K$ which is (AS).
This result is the same as the analyzed in \cite{EsEmEs}.
\item \textbf{Fractional logistic equation with harvest:}
This is equation (\ref{problemappal}) with $a=0$, $b=-\tfrac{r}{K}$ and $c=r-E$. Assuming that the intrinsic growth rate is bigger than the harvest coefficient $(r>E)$,  $a=0$, $b<0$ and $c>0$, then the equilibrium points are $x_1=0$ which is (U) and $x_2=K\left(1-\tfrac{E}{r}\right)$ which is (AS). In the case where the intrinsic growth rate is lower than the harvest coefficient  $(r<E)$, $a=0$, $b<0$ and $c<0$ is obtained, resulting $x_1=0$ (AS) and $x_2=K\left(1-\tfrac{E}{r}\right)$ (U).
\item \textbf{Fractional logistic equation with Allee effect:}
This is equation (\ref{problemappal}) with $a=-\tfrac{r}{K}$, $b=\left(\tfrac{m}{K}+1\right)r$ and $c=-rm$.
Then $a<0$, $b>0$ and $c<0$, therefore the equilibrium points are $x_1=0$ which is (AS), $x_2=m$ which is (U), because $b>\sqrt{b^2-4ac} $, and $x_3=K$ which is (AS), because $b>-\sqrt{b^2-4ac}$. This result is the same as the analyzed in \cite{SyMaSh}.
\item \textbf{Fractional logistic equation with Allee effect with harvest:}
This is equation (\ref{problemappal}) with $a=-\tfrac{r}{K}$, $b=\left(\tfrac{m}{K}+1\right)r$ and $c=-rm-E$.
Then $a<0$, $b>0$ and $c<0$. The equilibrium points are classified according to the value of $E$. If $E<\frac{r}{4K}(K-m)^2$ then they are $x_1=0$ which is (AS), $x_2=\frac{1}{2}\left(m+K+\sqrt{(K-m)^2-4\tfrac{KE}{r}}\right)$ which is (U), because $b>\sqrt{b^2-4ac} $, and $x_3=\frac{1}{2}\left(m+K-\sqrt{(K-m)^2-4\tfrac{KE}{r}}\right)$ which is (AS), because $b>-\sqrt{b^2-4ac}$, while if $E \geq \frac{r}{4K}(K-m)^2$  the only equilibrium point is $x_1=0$ which is (AS), leading to the extinction of the species.
\end{itemize}

\subsection{Numerical results}
\,\\
To perform numerical implementations a predictor-corrector method is used. The  fractional forward Euler method is utilized to get $u_{n+1}^{P}$ (predictor), and then the fractional trapezoidal rule is used to get $u_{n + 1}$ (corrector), which leads to the fractional Adams method, \cite{BaDiScTru, LiZe}.
\[ \left\{ \begin{array}{ll}
u_{n+1}^{P}&=\displaystyle \sum_{j=0}^{m-1}\frac{t_{n+1}^{j}}{j!}u_{0}^{j}+\sum_{j=0}^{n}b_{j,n+1}f(t_j,u_j),\\
u_{n+1}&=\displaystyle \sum_{j=0}^{m-1}\frac{t_{n+1}^{j}}{j!}u_{0}^{j}+\sum_{j=0}^{n}a_{j,n+1}f(t_j,u_j)+a_{n+1,n+1}f(t_{n+1},u_{n+1}^{P}).\\
\end{array}\right.\]

Different graphics of $ x(t) $ are presented taking into account different variations of the parameters $r,\,K,\,m,\,E,\,x_0,\,\alpha$ and $T$.

\begin{itemize}
\item In the following graphics, the fractional logistic equation with harvest considered is
\[\,_0^{C} D_t^{0.5}\left[x\right](t)= 0.5\, x(t) \left(1-\frac{x(t)}{10}\right)-E \,x(t)\]

where $r=0.5$ is the same as in \cite{EsEmEs} but considering a carrying capacity of $K=10$ and $\alpha=0.5$. Each graphic represents the approximate solution of the fractional equation for a given harvest value $E$, where the initial values $x(0)=x_0$ are varied. The chosen final time is $T=500$ since in general the solutions of the fractional logistic equation take a long time to reach the equilibrium.

\begin{figure}[H]
    \centering
    \includegraphics[width=.4\textwidth]{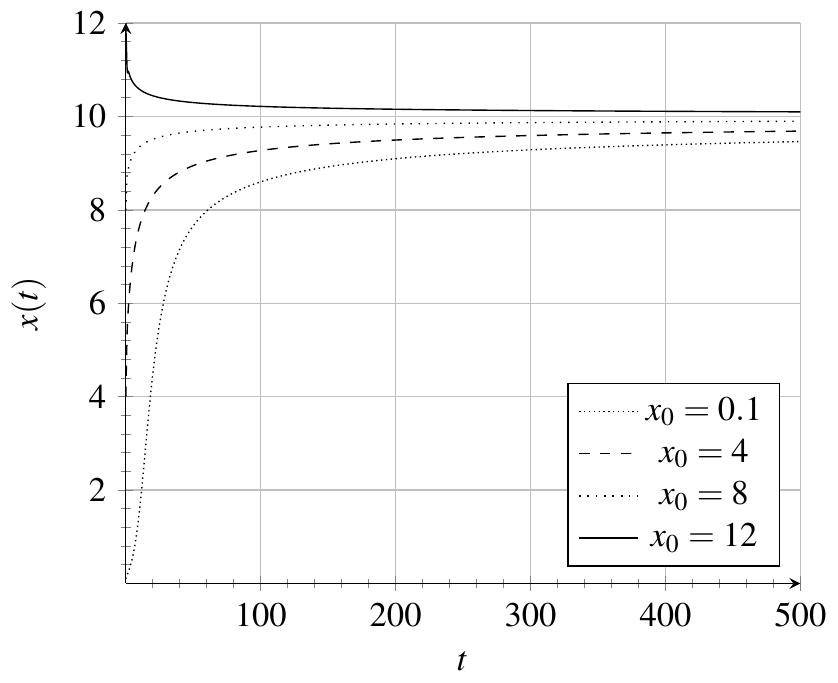} \includegraphics[width=.4\textwidth]{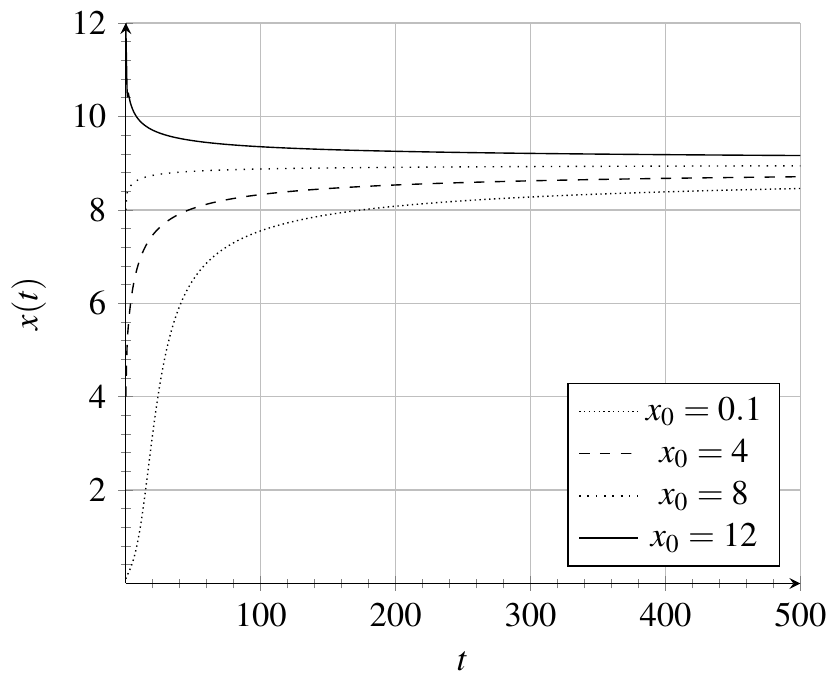}
    
    \includegraphics[width=.4\textwidth]{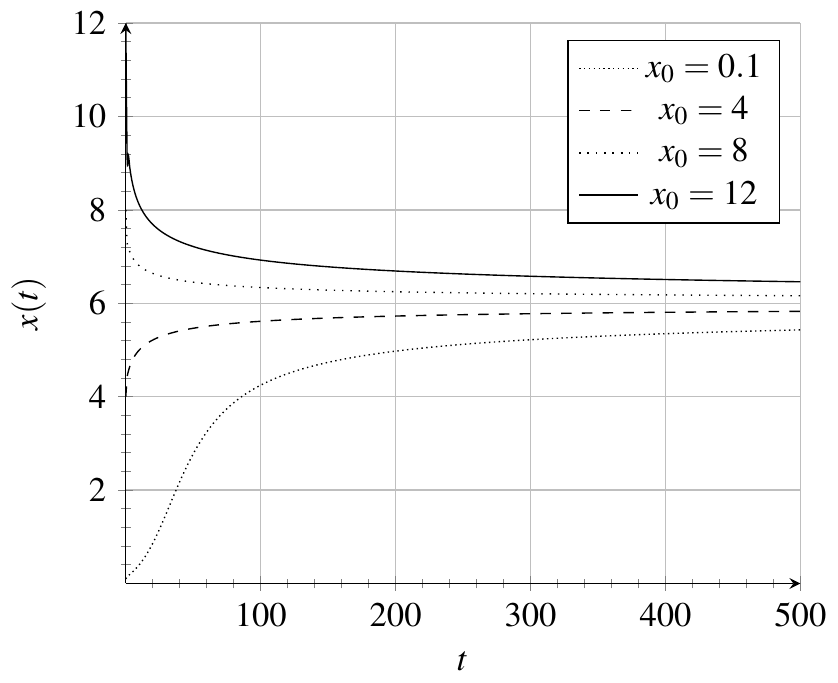} \includegraphics[width=.4\textwidth]{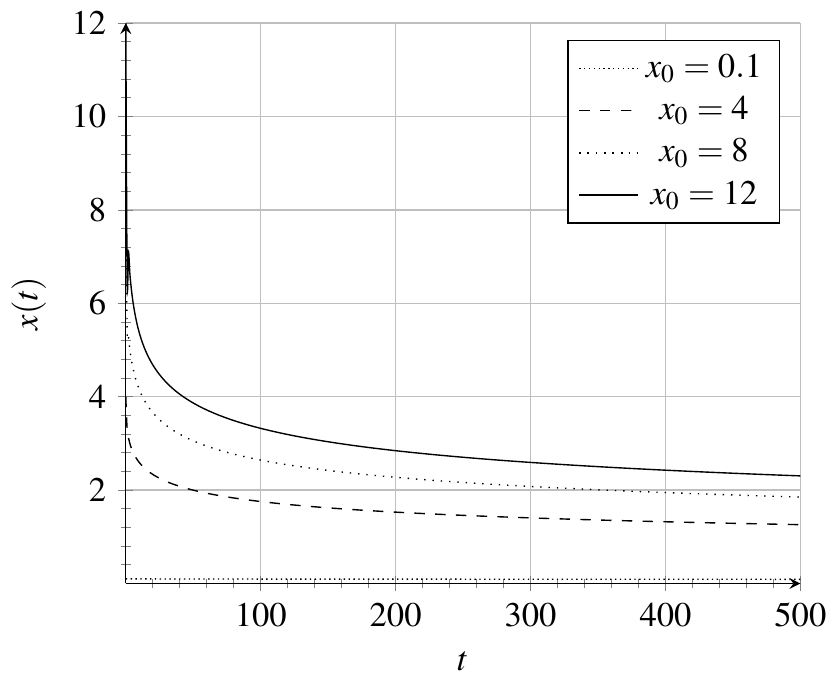}
    \caption{Fractional logistical equation solutions with $E=0, \, E=0.05, \, E=0.2$ and $E=0.5$}
    \label{Figlogcosecha1}
\end{figure}

It can be noticed in Figure \ref{Figlogcosecha1} that if the harvest of the species is zero, the solutions are stabilized in $K=10$, analogously to the behavior seen in \cite{EsEmEs} in which the solutions are also stabilized in the carrying capacity. As the harvest increases, there is stability at intermediate points until the species are extinguished when there is over-exploitation as shown in the last graphic. This is similar to the integer case, the difference is the way in which these solutions reach that equilibrium. 

\item  In the following graphics, the fractional logistic equation with harvest considered is
\[\,_0^{C} D_t^{\alpha}\left[x\right](t)= 0.5\, x(t) \left(1-\frac{x(t)}{10}\right)-0.2 \,x(t)\]

where again $r=0.5$, $K=10$ and $E=0.2$. In this case, each graphic represents the fractional equation approximate solution for a given initial value $x(0)=x_0$, where the values of $\alpha$ are varied. Again, the chosen final time is $T=500$.

\begin{figure}[H]
    \centering
    \includegraphics[width=.4\textwidth]{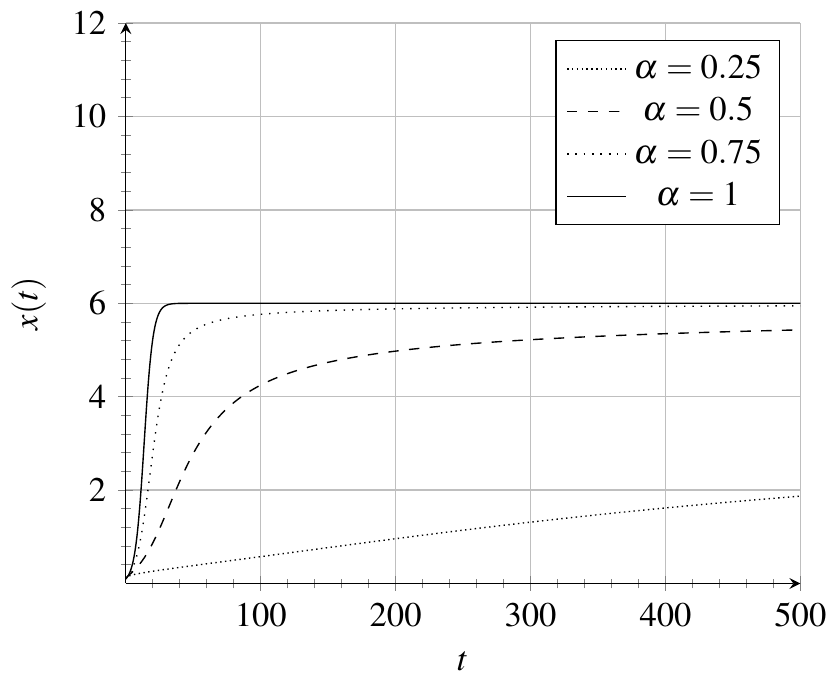} \includegraphics[width=.4\textwidth]{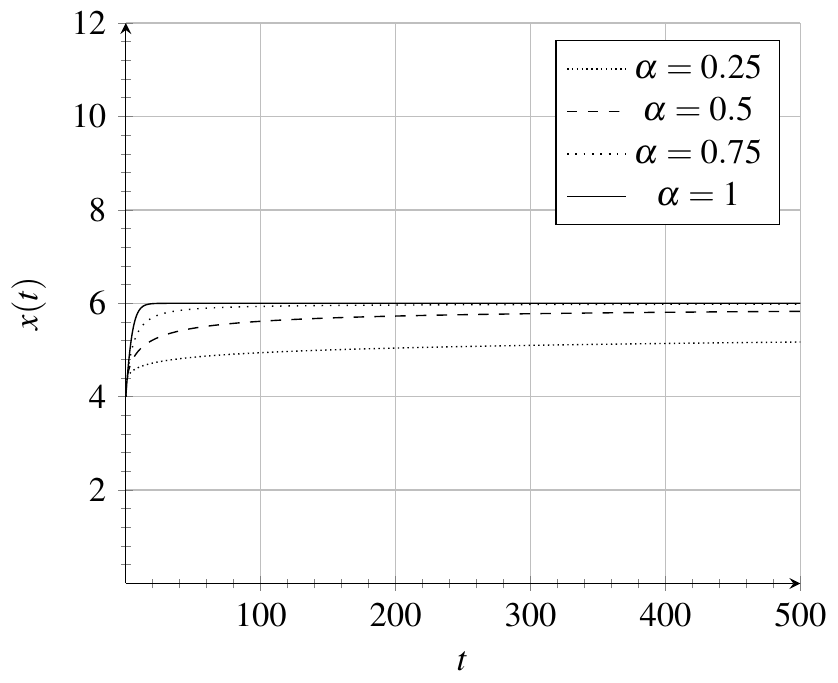}
    
    \includegraphics[width=.4\textwidth]{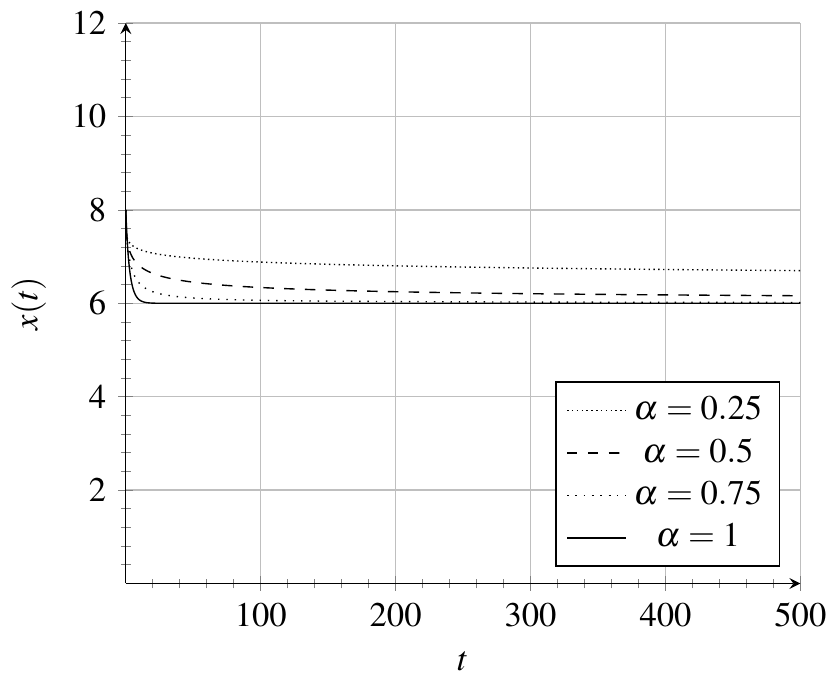} \includegraphics[width=.4\textwidth]{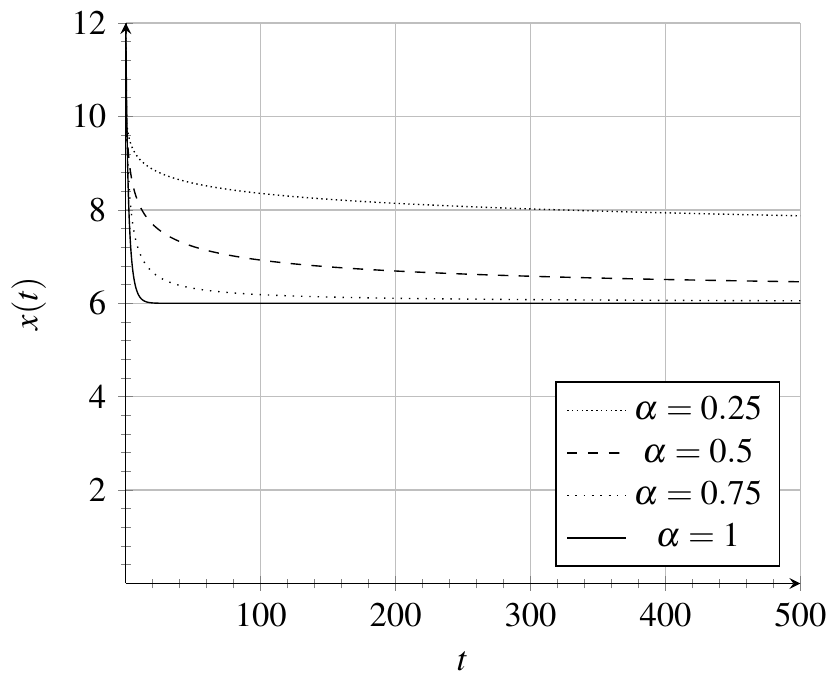}
    \caption{Fractional logistic equation with harvest solutions with $x_0=0.1, \, x_0=4, \, x_0=8$ and $x_0=12$}
    \label{Figlogcosecha2}
\end{figure}
It can be noticed in Figure \ref{Figlogcosecha2} that for different values of $\alpha$ the same equilibrium is achieved, but in a different way. The shapes of each solution $x(t)$ are different and the times in which each one reaches the equilibrium value are different as well, that is when the value of $\alpha$ decreases, the solutions reach equilibrium more slowly.

\item In the following graphics, the fractional logistic equation with Allee effect with harvest considered is
\[\,_0^{C} D_t^{0.5}\left[x\right](t)= 0.5\, x(t) \left(1-\frac{x(t)}{10}\right)\left(x(t)-1\right)-Ex(t)\]
where $r=0.5$, $K=10$, $m=1$, are the same as the parameters in \cite{SyMaSh} and $\alpha=0.5$. In this case, each graphic represents the fractional equation approximate solution for a given harvest value $E$, where the initial values $x(0)=x_0$ are varied. The chosen final time is $T=25$, because the solutions of the equation with Allee effect reach the equilibrium faster than the previous cases with no Allee effect.
\begin{figure}[H]
    \centering
    \includegraphics[width=.4\textwidth]{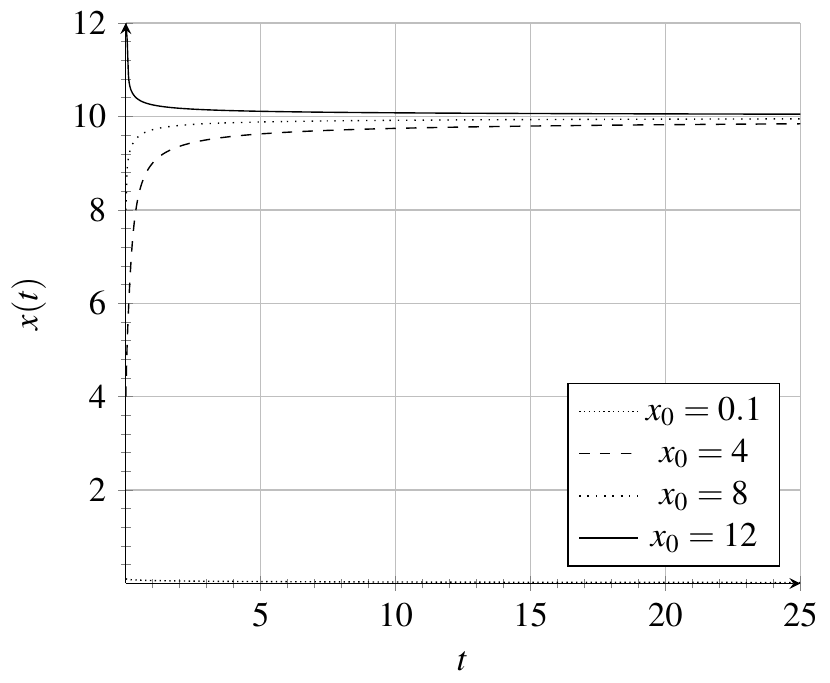} \includegraphics[width=.4\textwidth]{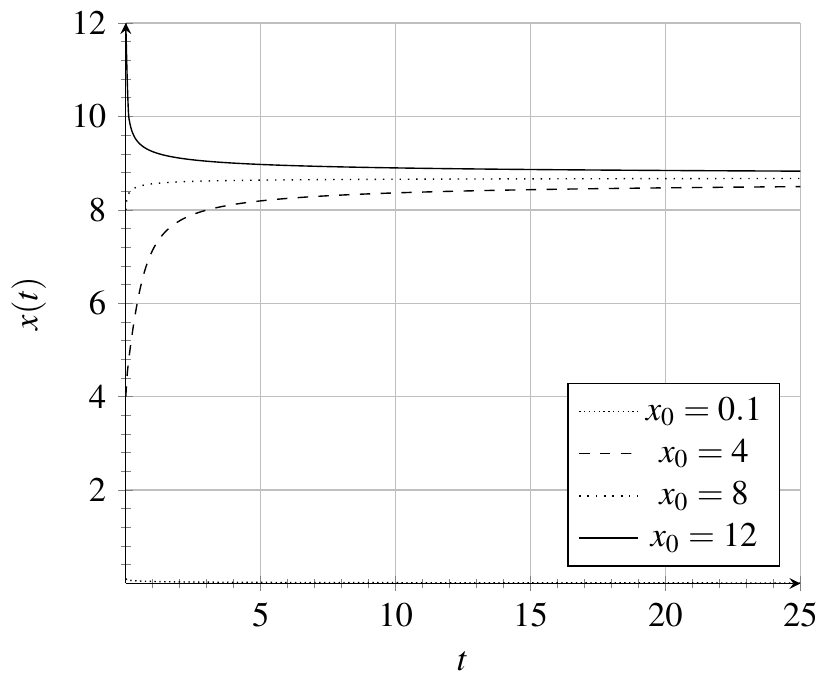}
    
    \includegraphics[width=.4\textwidth]{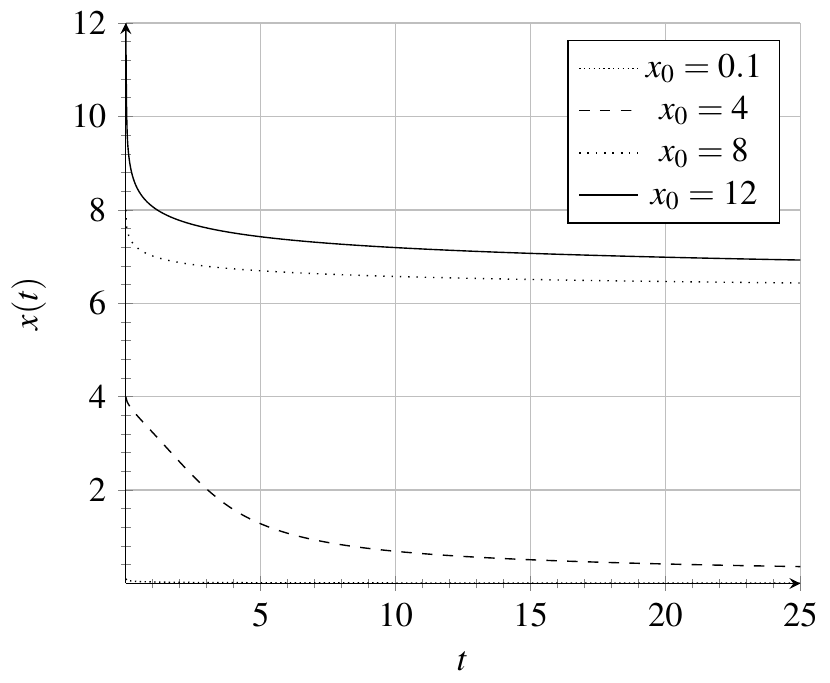} \includegraphics[width=.4\textwidth]{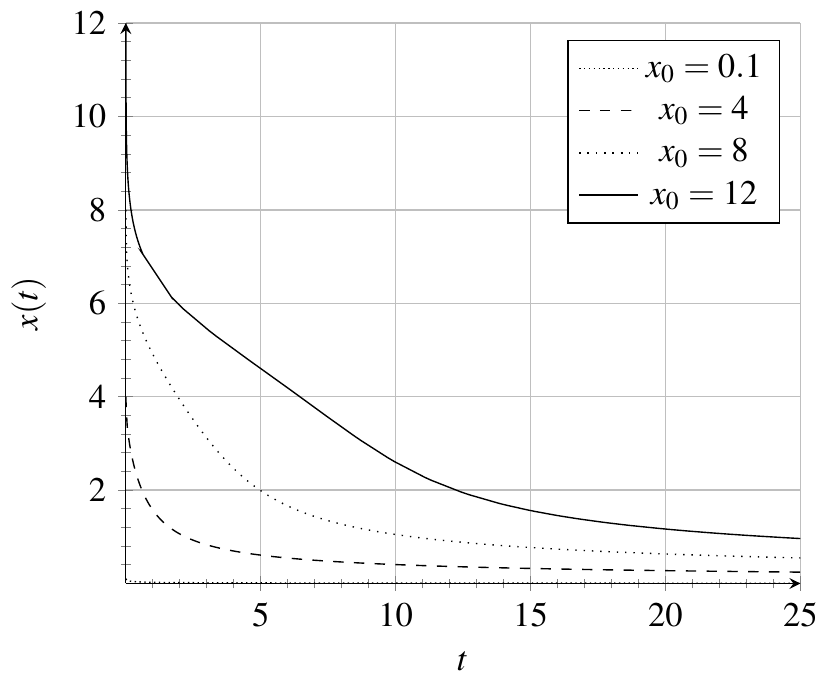}
    \caption{Fractional logistic equation with Allee effect with harvest solutions with ${E=0}, \, {E=0.5}, \, {E=1}$ and ${E=1.5}$}
    \label{Figalleecosecha1}
\end{figure}
It can be noticed in Figure \ref{Figalleecosecha1} that if the harvest is zero, the solutions stabilize in the carrying capacity $K=10$, analogously to \cite{SyMaSh}. As the harvest increases, there is stability at intermediate points until the species are extinguished when the harvest value $E=1.5>1.025=\frac{r}{4K}(K-m)^2$. This is similar to the integer case, the difference is the way in which these solutions reach this equilibrium, which is slower as the value of $\alpha$ decreases.

\item In the following graphics, the same fractional logistic equation with Allee effect with harvest considered is
\[\,_0^{C} D_t^{\alpha}\left[x\right](t)= 0.5\, x(t) \left(1-\frac{x(t)}{10}\right)\left(x(t)-1\right)-0.2\, x(t)\]

where again $r=0.5$, $K=10$, $m=1$ and $E=0.2$. In this case, each graphic represents the fractional equation approximate solution for a given initial value $x(0)=x_0$, where the values of $\alpha$ are varied. Again the chosen final time is $T=25$.
\begin{figure}[H]
    \centering
    \includegraphics[width=.4\textwidth]{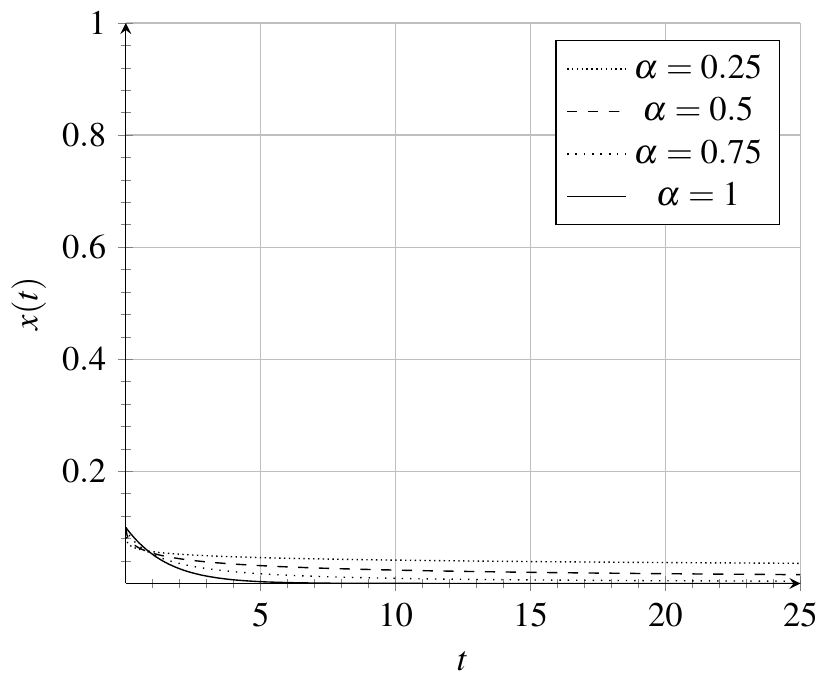} \includegraphics[width=.4\textwidth]{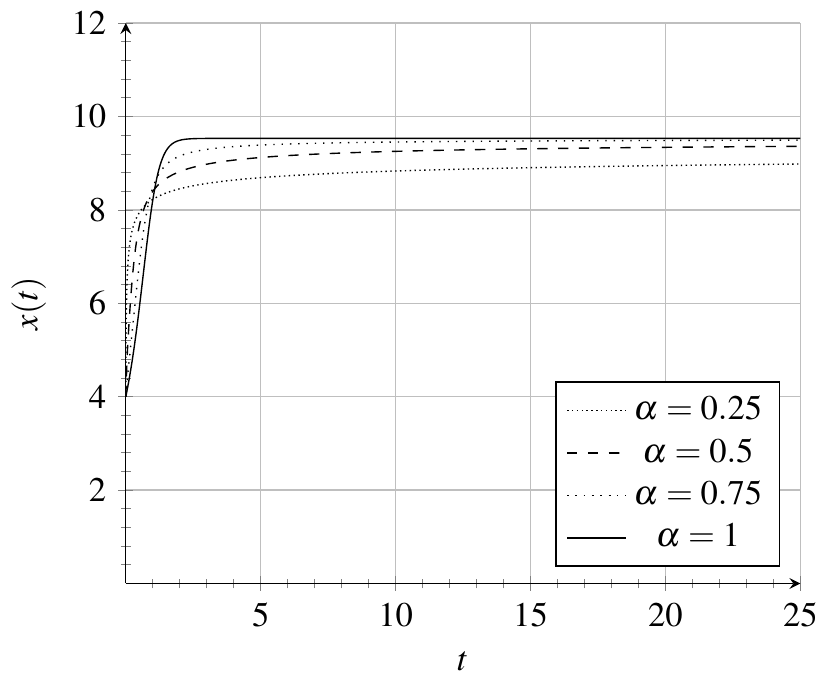}
    
    \includegraphics[width=.4\textwidth]{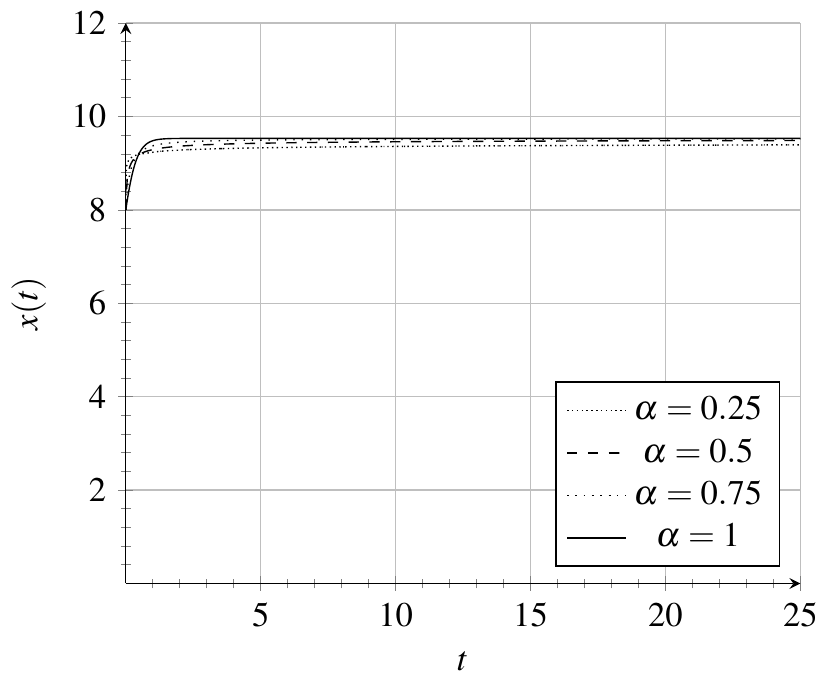} \includegraphics[width=.4\textwidth]{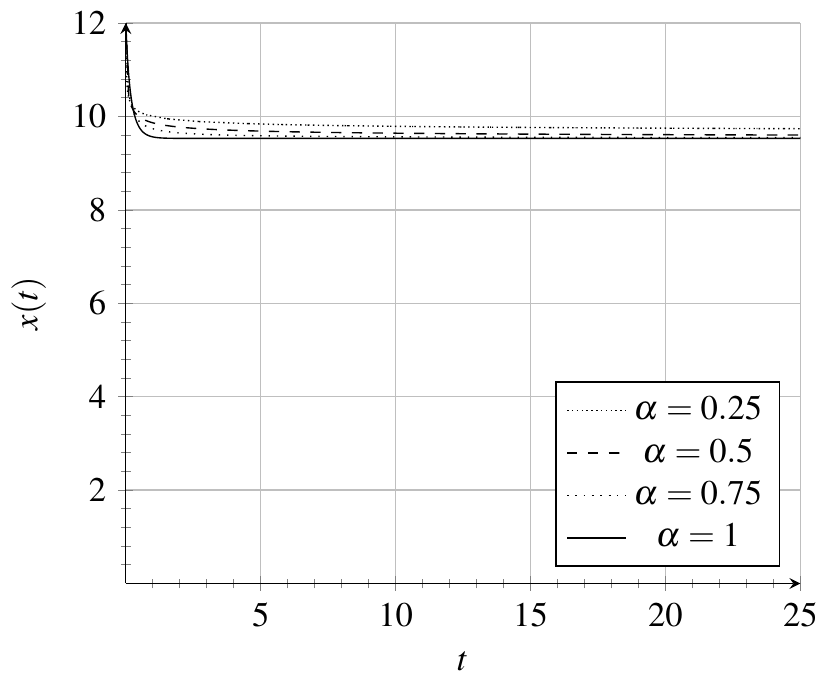}
    \caption{Fractional logistic equation with Allee effect with harvest solutions with $x_0=0.1, \, x_0=4, \, x_0=8$ and $x_0=12$}
    \label{Figalleecosecha2}
\end{figure}

It can be noticed in Figure \ref{Figalleecosecha2} hat for different values of $\alpha$  the same equilibrium is reached, but in a different way. The shapes of each solution $x(t)$ are different and also the times in which each one reaches the equilibrium value are different.
\end{itemize}

\section{Conclusions}
\label{sec:conclusions}
In this paper, an analysis of a fractional cubic equation, which is a generalization of different versions of fractional logistic equations, has been presented. In the first part of the work it was demonstrated the existence and uniqueness of the fractional cubic problem subject to initial values. In the second part a stability analysis was performed and it could be seen that this analysis is similar to the analysis of the fractional logistic equations. Finally, the fractional Adams method was numerically implemented, being able to extend results already obtained in order to compare the different fractional logistic models with harvest.

\section*{Acknowledgments}
 This work was partially supported by Universidad Nacional de Rosario through the projects ING568 ``Problemas de Control \'Optimo Fraccionario''. The first author was also supported by CONICET through a PhD fellowship.

\bibliographystyle{acm}
\bibliography{biblio} 

\section*{Appendix}
Consider $0<\alpha \leq 1$ and $a,b,c \in \mathbb{R}$ for the following problem,

\begin{equation}
 \left\{ \begin{array}{l}
             \,_0^{C} D_t^{\alpha}\left[x\right](t)=f(x(t))=ax^3(t)+bx^2(t)+cx(t), \\
             \\ x(0)=x_0. \\
             \end{array}
   \right.
	\label{problemappal20}
\end{equation}

To find the equilibrium points of the equation (\ref{problemappal20}), $\,_0^{C} D_t^{\alpha}x(t)=0$ is setted, therefore

\[ax^3(t)+bx^2(t)+cx(t)=0.\]

$\diamond$ If $a\neq 0$, three equilibrium points are obtained:

\[x_1=0 \, \, \text{ y } \, \, x_{2,3}=\frac{-b\pm \sqrt{b^2-4ac}}{2a}.\]

 Analyzing the sign of $\frac{\partial f}{\partial x}(x_{i}), \ i=1,2,3$, the stability of the equilibrium points will be obtained.
 
\begin{itemize}
\item \underline{Equilibrium point $x_1=0$}: as $\frac{\partial f}{\partial x}(x_{1})=c$, it is concluded that
\begin{itemize}
\item If $c > 0$, results $arg(\frac{\partial f}{\partial x}(x_{1}))=0\leq \tfrac{\alpha \pi}{2}$ then $x_1$ is (U).
\item If $c<0$, results $arg(\frac{\partial f}{\partial x}(x_{1}))=\pi>\tfrac{\alpha \pi}{2}$ then $x_1$ is (AS).
\item If $c=0$, there is no conclusion because $\frac{\partial f}{\partial x}(x_{1})=0$.
\end{itemize}

\item \underline{Equilibrium point $x_2=\frac{-b+\sqrt{b^2-4ac}}{2a}$}: as $\frac{\partial f}{\partial x}(x_{2})=\frac{b^2-b\sqrt{b^2-4ac}-4ac}{2a}$, it is concluded that

\begin{itemize}
\item If $a>0$ and $b<\sqrt{b^2-4ac}$, results $arg(\frac{\partial f}{\partial x}(x_{2}))=0<\tfrac{\alpha \pi}{2}$ then $x_2$ is (U).
\item If $a<0$ and $b>\sqrt{b^2-4ac}$, results $arg(\frac{\partial f}{\partial x}(x_{2}))=0<\tfrac{\alpha \pi}{2}$ then $x_2$ is (U).
\item If $a>0$ and $b>\sqrt{b^2-4ac}$, results $arg(\frac{\partial f}{\partial x}(x_{2}))=\pi>\tfrac{\alpha \pi}{2}$ then $x_2$ is (AS).
\item If $a<0$ and $b<\sqrt{b^2-4ac}$, results $arg(\frac{\partial f}{\partial x}(x_{2}))=\pi>\tfrac{\alpha \pi}{2}$ then $x_2$ is (AS).
\item If $ b=\sqrt{b^2-4ac}$ or $b^2-4ac< 0$, then $x_2=x_1=0$, that was already analyzed. 
\item If $b^2-4ac= 0$, there is no conclusion because $\frac{\partial f}{\partial x}(x_{2})=0$. 
\end{itemize}

\item \underline{Equilibrium point $x_3=\frac{-b-\sqrt{b^2-4ac}}{2a}$}: as $\frac{\partial f}{\partial x}(x_{3})=\frac{b^2+b\sqrt{b^2-4ac}-4ac}{2a}$, it is concluded that

\begin{itemize}
\item If $a>0$ and $b<-\sqrt{b^2-4ac}$, results $arg(\frac{\partial f}{\partial x}(x_{3}))=\pi>\tfrac{\alpha \pi}{2}$ then $x_3$ is (AS).
\item If $a<0$ and $b>-\sqrt{b^2-4ac}$, results $arg(\frac{\partial f}{\partial x}(x_{3}))=\pi>\tfrac{\alpha \pi}{2}$ then $x_3$ is (AS).
\item If $a>0$ and $b>-\sqrt{b^2-4ac}$, results $arg(\frac{\partial f}{\partial x}(x_{3}))=0<\tfrac{\alpha \pi}{2}$ then $x_3$ is (U).
\item If $a<0$ and $b<-\sqrt{b^2-4ac}$, results $arg(\frac{\partial f}{\partial x}(x_{3}))=0<\tfrac{\alpha \pi}{2}$ then $x_3$ is (U).
\item If $ b=-\sqrt{b^2-4ac}$ or $b^2-4ac<0$, then $x_3=x_1=0$, that was already analyzed.
\item If $b^2-4ac=0$, there is no conclution because $\frac{\partial f}{\partial x}(x_{3})=0$.
\end{itemize}

\end{itemize}
$\diamond$ If $a=0$, the equilibrium points are
\[x_1=0 \,\, \text{ y } \,\, x_{2}=-\frac{c}{b},\]
assuming that $b\neq 0$ and $c \neq 0$ since if this does not happen $x_1=0$ would be the only equilibrium point.\\
\noindent Analyzing the sign of $\frac{\partial f}{\partial x}(x_{i}), \ i=1,2$, the stability of the equilibrium points will be obtained.

\begin{itemize}
\item \underline{Equilibrium point $x_1=0$}: Analogous to the previous case.
\item \underline{Equilibrium point $x_2=-\frac{c}{b}$}: as $\frac{\partial f}{\partial x}(x_{2})=-c$, it is concluded that
\begin{itemize}
\item If $c>0$, results $arg(\frac{\partial f}{\partial x}(x_{2}))=\pi>\tfrac{\alpha \pi}{2}$ then $x_2$ is (AS).
\item  If $c<0$, results $arg(\frac{\partial f}{\partial x}(x_{2}))=0\leq\tfrac{\alpha \pi}{2}$ then $x_2$ is (U).
\end{itemize}
\end{itemize}
\end{document}